\documentclass[dvipsnames,letterpaper, 10 pt, journal]{ieeetran}

\IEEEoverridecommandlockouts	% This command is only
\usepackage{cite}
\usepackage{amsmath,amssymb,amsfonts,amsthm}
\usepackage{algorithmic}
\usepackage{textcomp}
\usepackage{graphicx,color}
\usepackage{mathrsfs}
\usepackage{tikz}
\usepackage[vlined,ruled]{algorithm2e}
\usepackage{subfigure}
\usepackage{url}
\usepackage{color}
\usepackage{dsfont}
\usepackage{bbm}
\usepackage{booktabs}
\usepackage{array}
\usepackage{yfonts}
\usepackage[normalem]{ulem}
\usepackage{arydshln,leftidx,mathtools}
\usepackage{multicol}
\usepackage{amsmath}
\usepackage{stfloats}
\usepackage{comment}
\usepackage{subfigure}

\definecolor{myBlue}{RGB}{49,130,189}

\newtheorem{theorem}{Theorem}[section]

\newtheorem{example}{Example}
\newtheorem{appxlem}{Lemma}[subsection]

\newcommand{\map}[3]{#1: #2 \rightarrow #3}

\newcommand{\blkdiag}{\mathrm{blkdiag}}

\usepackage{tabularx}
\usepackage{multirow}
\usepackage{makecell}

\newcommand\aamsout{\bgroup\markoverwith{\textcolor{violet}{\rule[0.5ex]{2pt}{1pt}}}\ULon}

\newcommand{\Rank}{\operatorname{Rank}}

\newcommand{\col}{\operatorname{Col}}

\newcommand{\real}{\mathbb{R}}

\newcommand{\T}{\mathsf{T}} %or \top or \intercal

\newcommand{\mc}{\mathcal}

\newcommand{\vv}{\text{v}}
\newcommand{\h}{\text{h}}
\newcommand{\comp}{\text{c}}

\newcommand{\vectorize}[1]{\text{vec}\!\left(#1\right)}
\newcommand{\vectorizee}[1]{\text{vec}\left(#1\right)}
\newcommand{\invvectorize}[1]{\text{vec}^{-1}\!\left(#1\right)}

\newcommand{\Range}[1]{\text{Im}\left(#1\right)}
\newcommand{\Null}[1]{\text{Ker}\left(#1\right)}

\newcommand{\1}{\mathds{1} }

\DeclareSymbolFont{bbold}{U}{bbold}{m}{n}
\DeclareSymbolFontAlphabet{\mathbbold}{bbold}

\newcommand\oprocendsymbol{\hbox{$\square$}}
\newcommand\oprocend{\relax\ifmmode\else\unskip\hfill\fi\oprocendsymbol}

\renewcommand{\theenumi}{(\roman{enumi})}

\newcommand*{\QEDA}{\hfill\ensuremath{\blacksquare}}%

\DeclareMathOperator{\diag}{diag}

\graphicspath{{figs/}}

\makeatletter
\let\NAT@parse\undefined
\makeatother
\usepackage[colorlinks,urlcolor=blue, linkcolor=blue, citecolor=blue]{hyperref}

\begin{document}

\title{\LARGE \bf Online Optimization with Unknown Time-varying Parameters}
%\title{\LARGE \bf Time-varying Polynomial Optimization with Unknown
%  Parameters using Gradient Measurements}

% \title{\bf Performance Limitations of Perception-Based Estimation in
%   Non-Nominal Environments}

\author{Shivanshu~Tripathi, Abed~AlRahman~Al~Makdah, and
  Fabio~Pasqualetti \thanks{This material is based upon work supported
    in part by awards ARO-W911NF-20-2-0267 and
    AFOSR-FA9550-20-1-0140. S.~Tripathi is with the Department of
    Electrical and Computer Engineering, and F. Pasqualetti are with
    the Department of Mechanical Engineering at the University of
    California, Riverside,
    \href{mailto:strip008@ucr.edu}{\{\texttt{strip008}},\href{mailto:fabiopas@ucr.edu}{\texttt{fabiopas\}@ucr.edu}}. A. A. Al
    Makdah is with the School of Electrical, Computer and Energy
    Engineering, Arizona State University,
    \href{mailto:aalmakda@asu.edu}{\texttt{aalmakda@asu.edu}}.}}

\maketitle

\begin{abstract}
  In this paper, we study optimization problems where the cost
  function contains time-varying parameters that are unmeasurable and
  evolve according to linear, yet unknown, dynamics. We propose a
  solution that leverages control theoretic tools to identify the
  dynamics of the parameters, predict their evolution, and ultimately
  compute a solution to the optimization problem. The identification
  of the dynamics of the time-varying parameters is done online using
  measurements of the gradient of the cost function. This system
  identification problem is not standard, since the output matrix is
  known and the dynamics of the parameters must be estimated in the
  original coordinates without similarity
  transformations. Interestingly, our analysis shows that, under mild
  conditions that we characterize, the identification of the
  parameters dynamics and, consequently, the computation of a
  time-varying solution to the optimization problem, requires only a
  finite number of measurements of the gradient of the cost
  function. We illustrate the effectiveness of our algorithm on a
  series of numerical examples. 
\end{abstract}

% Note that keywords are not normally used for peerreview papers.
%\begin{IEEEkeywords}
%
%\end{IEEEkeywords}

\section{Introduction}\label{sec: introduction}
Mathematical optimization is a fundamental tool used across diverse
applications, ranging from robotics \cite{PMW-MP-YH-AE-NM-ADP:23},
control \cite{DPB:18},~and machine learning \cite{CMB:06}. It offers a
structured framework for expressing objectives and defining
constraints, enabling a systematic approach to finding the most
favorable solutions in diverse contexts. Recently, there has been a
significant interest in time-varying optimization problems due to
their ability to capture the complexity of dynamic environments and
changing optimization objectives. In such problems, the~cost
parameters vary with time, and optimization algorithms must adapt in
real-time and incorporate new information to track the optimal
time-varying solution \cite{BTP:87,AS-ED-SP-GL-GBG:20}. In this work, we consider the following optimization problem with time-varying parameters
\begin{align} \label{eq: optimization prob_intro}
  \min_{x\in \mathbb{R}^n}~~f(x,z(t))=  g(x)^{\T} z(t), \quad t\geq 0 
\end{align}
where $\map{f}{\mathbb{R}^n \times \mathbb{R}^p}{\mathbb{R}}$ is the
cost function with time varying parameters $z(t)\in \mathbb{R}^p$, and
$\map{g}{\mathbb{R}^n}{\mathbb{R}^p}$ is a vector-valued function
whose entries depend on $x$. The parameters $z(t)$ are not measurable
and evolve with linear, yet unknown, dynamics. Our objective is to
develop an online algorithm to compute and track the optimal solution
of \eqref{eq: optimization prob_intro} using measurements of the
gradient of the cost function. In particular, we propose a solution
that leverages control theoretic tools to identify the exact dynamics
of the parameters $z(t)$, predict their evolution, and ultimately
compute the optimal solution. Interestingly, our analysis shows that
the identification of the cost function and, consequently, the
computation of the optimal time-varying solution to \eqref{eq:
  optimization prob_intro} requires only a finite number of
measurements of the gradient of the cost function.

\textbf{Related work.}  Several methods for solving time-varying
optimization problems have been studied in the literature. These
methods can be categorized into two classes, namely
\emph{unstructured} methods and \emph{structured}
methods. Unstructured methods are passive and do not leverage the
dynamics of the optimization problem: each time-instance of the
optimization problem is treated as a static problem and solved using
static optimization algorithms (e.g., gradient descent)
\cite{SSS:12,ECH-RMW:15,RD-ASB-RT-KR:19,SMF:20}. Since these methods
ignore the dynamic nature of the problem and make decisions only after
observing a time instance of the cost, convergence to the optimal
time-varying solution is never achieved and is only guaranteed to a
neigberhood of the optimal solution
\cite{ED-AS-SB-LM:20,AS-ED-SP-GL-GBG:20,PCV-SJ-FB:22,AD-VC-AG-GR-FB:23}. On
the counter part, structured algorithms leverage the dynamics of the
problem in order to achieve exact tracking of the time-varying optimal
solution \cite{AS-ED-SP-GL-GBG:20}. A particular class of algorithms
that fall under structured methods is the prediction-correction
algorithms. At each time-step, a prediction step is used to forecast
how the optimal solution evolves with time, followed by a correction
step to update the solution. The prediction-correction algorithm has
been studied in the literature in continuous-time
\cite{MB-CL-UH:04,SR-WR:16,MF-SP-VMP-AR:17} and discrete-time
\cite{AS-AM-AK-GL-AR:16,AS-ED:17,AS:18}. In the recent work
\cite{NB-RC-SZ:24}, the authors study quadratic optimization problems
with time-varying linear term, where they leveraged control theoretic
tools to develop a structured online algorithm that achieves
zero-tracking error asymptotically. In \cite{UC-NB-RC-SZ:24}, the
authors extended the work of \cite{NB-RC-SZ:24} by studying quadratic
optimization problems with inequality constraints where the linear
term in the cost and the inequality values are time-varying. In this
work, we focus on optimization problems as in \eqref{eq: optimization
  prob_intro}, where we propose a solution that leverages control
theoretic tools to identify the exact dynamics of the coefficients,
$z(t)$, predict their evolution, and ultimately compute the optimal
solution. However, differently from \cite{NB-RC-SZ:24,UC-NB-RC-SZ:24},
we have all the cost parameters $z(t)$ in \eqref{eq: optimization
  prob_intro} varying with time as opposed to only varying the linear
term. Recent works \cite{GB-BVS:24,GB-BVS:24:1} have explored
time-varying cost functions where temporal variability arises from the
time-varying parameters. These works consider a general cost function
with unknown parameters that vary linearly, assuming either the
parameters are known or the dynamics of these parameters are known, or
both the parameters and their dynamics are known. In contrast, our
work addresses a more general setting where neither the parameters nor
their dynamics are known. However, we do restrict our analysis to cost
functions expressed as \eqref{eq: optimization prob_intro}.

\noindent
\textbf{Contributions.} The main contributions of this paper are as
follows. First, we design an algorithm to solve the time-varying
optimization problem \eqref{eq: optimization prob_intro}, which requires the
identification of the dynamics of the parameters in the cost using
measurements of the gradient of the cost function. This system
identification problem is not standard, since the output matrix is
known and the dynamics of the parameters must be estimated in the
original coordinates without similarity transformations. Second, we
provide a set of conditions on the parameters dynamics and gradient
measurements that guarantee the solvability of the parameters
identification problem using only a finite number of gradient
measurements. Third and finally, we illustrate the effectiveness of
our algorithm on a series of examples.

\noindent
\textbf{Notation.} The $n\!\times\! n$ identity matrix is denoted by
$I_n$. The~Kronecker product is denoted by
$\otimes$. The left (right) pseudo inverse of a tall (fat) matrix $A$
is denoted by $A^{\dagger}$. The column-space of a matrix $A$ is
denoted by $\col{\left(A\right)}$. The dimension of~a subspace
$\mc{V}$ is denoted by $\dim{\left(\mc{V}\right)}$. The range-space
and the null-space of a matrix $A$ are denoted by $\Range{A}$ and
$\Null{A}$, respectively. For a matrix
$A \!\in\! \mathbb{R}^{n \times m}$, let
$\vectorize{A} \!\in\! \mathbb{R}^{nm}$ denotes the vectorization of $A$,
which is obtained by stacking the columns of $A$ on top of one
another, and $\invvectorize{\cdot}$ denotes the inverse vectorization operator, i.e., $\invvectorize {\vectorize A}=A$. The complement~set and the
 cardinality of a set $\mc{S}$ are denoted by $\mc{S}^{\comp}$ and
 $|\mc{S}|$, respectively.

 \section{Problem formulation}\label{sec: formulation}
Consider the time-varying optimization problem \eqref{eq: optimization
  prob_intro}. Let the parameters $z(t)$ evolve according to
the following discrete-time, linear, time-invariant dynamics: 
\begin{align}\label{eq: parameter dynamics block}
  z(t+1)  = A z(t),
\end{align}
where $A\in \mathbb{R}^{p\times p}$. We assume that the matrix $A$ is
unknown, and that the parameters $z$ are not measurable. However, we
assume that, at each time, an Oracle can provide the value of the
gradient of the cost function $f$ in \eqref{eq: optimization
  prob_intro} with respect to $x$ evaluated at a desired point
$x \in \mathbb{R}^n$. Notice that
\begin{align}\label{eq: output model}
  \underbrace{\frac{\partial f}{\partial x}\bigg|_{x}}_{y(x,t)} = 
  \underbrace{\frac{\partial g^{\T}}{\partial x}\bigg|_{x}}_{C(x)}
  z(t) .
\end{align}
In other words, the gradient of the cost function evaluated at a
desired point $x$ is a linear function of the parameters $z$, which
can be used as a linear output to simultaneously identify the
parameters dynamics, predict their trajectory, and ultimately solve
the time-varying optimization problem \eqref{eq: optimization
  prob_intro}. It should be noticed that the identification of the
parameters dynamics is a nonstandard system identification problem,
because the output matrix $C(x)$ is known yet dependent on the
optimization variable $x$, and because the parameters dynamics
\eqref{eq: parameter dynamics block} need to be identified in their
original coordinates rather than up to an arbitrary similarity
transformation. In fact, a similarity transformation of \eqref{eq:
  parameter dynamics block} would yield an incorrect set of parameters
and ultimately an incorrect solution to \eqref{eq: optimization
  prob_intro}. We make the following assumptions on the parameter
dynamics \eqref{eq: parameter dynamics block} and their initial state
$z(0)$ that hold throughout the letter:

\begin{itemize}
  
\item[(A1)] The matrix $A$ in \eqref{eq: parameter dynamics block} has distinct and non-zero eigenvalues.

\item[(A2)] The pair $(A,z(0))$ is controllable.
  
\end{itemize}
While Assumption (A1) is a technical condition that simplifies the
notation and analysis, Assumption (A2) is generically
satisfied and necessary for the solvability of the identification and
estimation problems. Also, Assumption (A2) ensures that the initial condition $z(0)$ excites all the modes of the~system~\eqref{eq: parameter dynamics block}.

\begin{example}
Consider the time-varying optimization problem,
\begin{align*}
 \min_{x\in \mathbb{R}^2}~f(x,z(t)) =&z_1(t)x_1 + z_2(t)x_1x_2 + z_3(t)x_1^2 + z_4(t)x_2^2 \\
& + z_5(t)\cos(x_1) + z_6(t)\sin(x_2),  \quad t\geq0,
\end{align*}
where $\map{f}{\mathbb{R}^2 \times \mathbb{R}^6}{\mathbb{R}}$ is the cost function, $x=[x_1,x_2]^{\T}$ is the optimization variable, and $z(t)=[z_1(t), z_2(t), z_3(t), z_4(t), z_5(t), z_6(t)]^{\T}$ are the time-varying parameters of the cost. We can write $f(x,z(t))=g(x)^{\T} z(t)$ with $g(x)^{\T}= \begin{bmatrix}
 x_1 & x_1x_2 & x_1^2 & x_2^2 & \cos(x_1) & \sin(x_2)
\end{bmatrix}$. Let $z(t)$ evolve with time according to the linear dynamics in \eqref{eq: parameter dynamics block}, then, we can write the gradient measurement evaluated at $x$~as
\begin{align*}
   \underbrace{\frac{\partial f}{\partial x}\bigg|_{x}}_{y(x,t)} =
   \underbrace{
\begin{bmatrix}
 1 & x_2 & 2x_1 & 0 & -\sin(x_1) & 0\\
 0 & x_1 & 0       & 2x_2 & 0& \cos(x_2)
\end{bmatrix}}_{C(x)}
 z(t).
\end{align*}~\oprocend
\end{example}

\section{Estimating the parameters and their dynamics}\label{sec3}
In this section we detail our method to identify the unknown parameters and their
dynamics \eqref{eq: parameter dynamics block} to solve the
time-varying optimization problem \eqref{eq: optimization prob_intro}. Our
approach consists of three mains steps, namely, (i) the identification
of the parameter dynamics up to a similarity transformation using
subspace identification techniques, (ii) the identification of the
similarity transformation to recover the exact parameters dynamics using
the knowledge of the output matrix $C(x)$, and (iii) the prediction of
the parameters trajectory and the computation of a time-varying 
minimizer of \eqref{eq: optimization prob_intro}. For notational convenience, let $y(t)=y(x(t),t)$, where $x(t)$ is the
value at time $t$ at which the gradient $y(x,t)$ is evaluated. Let $X$
and $Y$ be the data collected from the Oracle over time, where
$N \in \mathbb{N}$ and
\begin{align}\label{eq: data}
  X=
  \begin{bmatrix}
    x(0)  & \cdots & x(N)
  \end{bmatrix}, \quad
                     Y=
                     \begin{bmatrix}
                       y(0) & \cdots & y(N)
                     \end{bmatrix}.
\end{align}
Further, let the data matrices $X$ and $Y$ be partitioned as
\begin{align}\label{eq: data partitioned}
 X=
\begin{bmatrix}
 X_0 & X_1
\end{bmatrix}, \quad
Y=
\begin{bmatrix}
 Y_0 & Y_1
\end{bmatrix},
\end{align}
where $X_0$ and $Y_0$ contain the first $N_0=2p-2$ samples of $X$
and $Y$, respectively, and $X_1$ and $Y_1$ the
remaining samples.
We make the following assumptions:
\begin{itemize}
\item[(A3)] The pair $(A, C(x_0))$ is observable, where $x_0$ is the
  point at which the gradient $y(x,t)$ is evaluated at time $t = 0$.
\item[(A4)] The values $x(t)$ remain constant up to time $N_0$, that
  is, for $0 \le t \le N_0$ and some vector $x_0 \in \real^n$,
  $x(t) = x_0$.
\end{itemize}
Assumption (A3) is generically
satisfied and necessary for the solvability of the identification and
estimation problems. Assumption (A4) implies that the output matrix $C(x)$ and the data
matrix $X_0$ are constant. This allows for the use of standard
subspace identification techniques and simplifies the derivations, and
it is compatible with settings where the Oracle can arbitrarily select
the gradient evaluation points. Additionally, since $x(t)$ remains
constant only for a finite time interval, this assumption induces only
a finite delay in the computation of the solution to the optimization
problem. Note that the time $N_0$ depends linearly on the system
dimension $p$. That is, more complex cost functions (i.e., more
parameters) induce a larger delay. The use of a non-constant matrix
$C(x)$ and data matrix $X$ remains the topic of 
current investigations.
We rewrite $Y_0$ in \eqref{eq: data partitioned} in the following
Hankel matrix form:
\begin{align}\label{eq: hankel matrix}
  Y_{\h}=
  \begin{bmatrix}
    y(0) & \cdots & y(p-1)\\
    \vdots & \ddots & \vdots \\
    y(p-1) & \cdots & y(2p-2)
  \end{bmatrix}.
\end{align}
The above Hankel matrix satisfies the following relation:
\begin{align}\label{eq: Z0}
  Y_{\h}= 
  \underbrace{\begin{bmatrix}
      C(x_0)\\
      C(x_0)A\\
      \vdots \\
      C(x_0) A^{p-1}
    \end{bmatrix}}_{O(x_0)}
  \underbrace{\begin{bmatrix}
      z(0) & Az(0) & \cdots & A^{p-1} z(0)
    \end{bmatrix}}_{Z(0)}.
\end{align}
Since $C(x_0)$ and $z(0)$ satisfy Assumptions (A2) and (A3), we have
$\Rank{\left(Y_\h \right)} = p$. 
To estimate
$A$, we make use of the shift structure of the observability matrix
$O(x_0)$ as in \cite{SYK:78}, i.e.,
\begin{align}\label{eq: A id}
  A=O_1^{\dagger}O_2,
\end{align}
where $O_1$ and $O_2$ are matrices containing the first $n(p-1)$ rows of $O(x_0)$ and the last
$n(p-1)$ rows of $O(x_0)$, respectively. We
now factor $Y_{\h}$ in \eqref{eq: hankel matrix} using the singular
value decomposition, $Y_{\h}=U\Sigma V^{\T}$, and construct the
observability matrix from the left singular vectors, $U$,
corresponding to the non-zero singular values. Notice that, because of
Assumptions (A2) and (A3), we have
$\col{\left(Y_{\h}\right)} = \col{\left(O(x_0)\right)}$. Hence, we
have $\overline{O}(x_0)\!=\!O(x_0)T^{-1}\!=\!\overline{U}$, where
$T\! \in\! \mathbb{R}^{p\times p}$ is a transformation matrix and
$\overline{U}$ contains the left singular vectors of $Y_{\h}$
corresponding to the non-zero singular values. Then,
\begin{align}\label{eq: A_bar id}
  \overline A = (\overline{O}_1)^{\dagger} \overline{O}_2,
\end{align}
where the matrices $\overline{O}_1$ and $\overline{O}_2$ contain the
first $n(p-1)$ rows of $\overline{O}(x_0)$ and the last $\left.n(p-1)\right.$ rows of $\overline{O}(x_0)$, respectively. By
construction, $\overline{A}$ is similar to $A$ in the parameters
dynamics, that is, $A = T^{-1}\overline{A} T$, for some invertible
matrix $T$. Next, we estimate the matrix $T$.

Using the notation in \eqref{eq: Z0}, let
\begin{align}\label{eq: Z_bar}
  \overline{Z}(0) = T Z(0) = 
  \begin{bmatrix}
    \overline{z}(0) & \cdots & \overline{z}(p-1)
  \end{bmatrix}
                               = {\overline{O}(x_0)}^{\dagger} Y_{\h},
\end{align}
and, to simplify the notation, let $C(t) = C(x(t))$. Using
\cite[Proposition 7.1.9]{DSB:09} and letting $\overline{z} = T z$, we have
\begin{align*}
  y(t) = \overline{C}(t) \overline{z}(t) &= C(t) T^{-1} \overline{z}(t)
  \\
  &=\left(\overline{z}(t)^{\T} \otimes C(t) \right)\vectorize{T^{-1}} .
\end{align*}
Consequently, we can write $Y_1$ in \eqref{eq: data partitioned} as
\begin{align}\label{eq: relation vec(T)}
  \underbrace{\begin{bmatrix}
      y(N_0+1)\\
      \vdots\\
      y(N)
    \end{bmatrix}}_{Y_{\vv}=\vectorizee{Y_1}}
  =
  \underbrace{\begin{bmatrix}
      {\overline{z}(N_0+1)}^{\T} \otimes C(N_0+1)\\
      \vdots \\
      {\overline{z}(N)}^{\T} \otimes C(N)
    \end{bmatrix}}_{M}
  \vectorize{T^{-1}}.
\end{align}
Note that the matrix $M$ in \eqref{eq: relation vec(T)} is known since
we can compute $\overline{z}(N_0+1), \cdots , \overline{z}(N)$ by
propagating $\overline{z}(p-1)$ in \eqref{eq: Z_bar} using the
matrix $\overline{A}$ in \eqref{eq: A_bar id}, and
$C(N_0+1), \cdots, C(N)$ are computed using $X_1$ in \eqref{eq: data
  partitioned}. Then, using \eqref{eq: relation vec(T)}, we can
compute $T$ as
\begin{align}\label{eq: transformation matrix}
  T= \left( \invvectorize{M^{\dagger} Y_{\vv}}   \right)^{-1}.
\end{align}
and, finally, the exact realization $A$ in \eqref{eq: parameter
  dynamics block} as
\begin{align*}
  A=T^{-1}\overline{A} T.
\end{align*}
This procedure relies on the matrix $M$ in \eqref{eq: relation vec(T)}
being full column rank, which we characterize in the next results.

\begin{theorem}{\bf \emph{(Necessary condition for the rank of $M$ in \eqref{eq: relation
        vec(T)})}}\label{thrm: rank M} The matrix $M$ in \eqref{eq: relation vec(T)} is full column-rank only if the matrix ${[C{(N_0 +1)}^{\T} \cdots {C(N)}^{\T}]}^{\T}$ is full rank.
\end{theorem}
\begin{IEEEproof}
Suppose $M$ in \eqref{eq: relation vec(T)} is full column-rank. Then,
\begin{align}\label{eq: kernel intersection eq}
  \bigcap_{i=N_0+1}^{N} \Null{ \overline{z}(i)^{\T}\otimes C(i)}=0.
\end{align}
Using Lemma \ref{lemma: 1}, we re-write \eqref{eq: kernel intersection eq} as
\begin{align}\label{eq: kernel intersection 1}
&\bigcap_{i=N_0+1}^{N}\! \Null{\overline{z}(i)^{\T}\!\otimes\! I_p} \cup \Null{I_p\! \otimes \!C(i)}\!=\!0.
\end{align}
For notational convenience, we use $\mc{A}(i)$ and $\mc{B}(i)$ to denote $\Null{\overline{z}(i)^{\T}\!\otimes\! I_p}$ and $\Null{I_p\!\otimes \!C(i)}$ for $i\!\in\!\{N_0\!+\!1,\!\cdots\! ,N\}$, respectively. Expanding the terms in \eqref{eq: kernel intersection 1}, we can write
\begin{align}\label{eq: kernel intersection expanded}
\begin{split}
\bigcup_{\mc{S}\subseteq \{N_0+1,\cdots, N\}} \left(\bigcap_{i\in \mc{S}} \mc{A}(i) \cap \bigcap_{j\in \mc{S}^{\comp}} \mc{B}(j) \right)=0.
\end{split}
\end{align}
Since the union of all the terms in \eqref{eq: kernel intersection expanded} is zero, then, each term of the union in \eqref{eq: kernel intersection expanded} is zero. Hence, we have
\begin{align}
\bigcap_{i=N_0+1}^{N}\mc{B}(i)&=\bigcap_{i=N_0+1}^{N} \Null{I_p\otimes C(i)}=0. \label{eq: cond 2}
\end{align}
Let $K_i$ be a matrix whose columns span $C(i)$. Then, using Lemma \ref{lemma: 2}, equation \eqref{eq: cond 2} can be equivalently written as
\begin{align}
\bigcap_{i=N_0+1}^{N}\col{\left(I_p\otimes K_i\right)}&=0. \label{eq: cond 2-2}
\end{align}
Equation \eqref{eq: cond 2-2} implies that
\begin{align*}
\begin{bmatrix}
 I_p \otimes C(N_0+1)\\ \vdots \\  I_p \otimes C(N)
\end{bmatrix}(I_p\otimes K)
= 0 \iff K=0.
\end{align*}
Equivalently, we have
\begin{align*}
\begin{bmatrix}
C(N_0+1)\\ \vdots \\  C(N)
\end{bmatrix}K
= 0 \iff K=0.
\end{align*}
Which implies that ${[C{(N_0 +1)}^{\T} \cdots {C(N)}^{\T}]}^{\T}$ is full rank.
\end{IEEEproof}
\begin{theorem}{\bf \emph{(Necessary and sufficient condition for the rank of $M$ in \eqref{eq: relation
        vec(T)})}}\label{thrm: sufficient cond}
  Let $\overline{A}$ in \eqref{eq: A_bar id} be diagonalizable and let
  $D$ be its diagonal form. Let
  $\left.\Lambda = D \otimes I_p\right.$ and
  $F(t) = {\1_p}^{\T} \otimes C(t)$, for $t\geq0$. Then, the matrix
  $M$ in \eqref{eq: relation vec(T)} is full column-rank if and only
  if the following matrix is full column-rank:
  \begin{align*}
    W = 
    \begin{bmatrix*}[c]
      F(N_0+1)\\
      F(N_0+2)\Lambda\\
      \vdots\\
      F(N) \Lambda^{N-N_0-1}
    \end{bmatrix*}
    .
  \end{align*}
\end{theorem}
\oprocend
\begin{IEEEproof}
From Assumption (A1), we have $A$ in \eqref{eq: parameter dynamics block} is diagonalizable. Hence, $\overline{A}$ in \eqref{eq: A_bar id} is diagonalizable. Let $D$ be the diagonal form of $\overline{A}$, i.e., $D=U^{-1}\overline{A}U$, where $U$ is a matrix whose columns are the eigenvectors of $\overline{A}$. Let $\xi(t)=U^{-1}\overline{z}(t)$ for $t\geq0$. Then, we can write $M$ in \eqref{eq: relation
        vec(T)} as
\begin{align}\label{eq: Q matrix}
\begin{split}
M&=
\begin{bmatrix}
      {\xi(N_0+1)}^{\T} U^{\T}\otimes C(N_0+1)\\
      \vdots \\
      {\xi(N)}^{\T}U^{\T} \otimes C(N)
    \end{bmatrix}\\
    &=
    \underbrace{\begin{bmatrix}
      {\xi(N_0+1)}^{\T}\otimes C(N_0+1)\\
      \vdots \\
      {\xi(N)}^{\T} \otimes C(N)
    \end{bmatrix}}_{Q}
    (U^{\T}\otimes I_p).
    \end{split}
\end{align}
It follows from \cite[Fact 7.4.24]{DSB:09} that the matrix $U^{\T}\otimes I_p$ is full-rank. Then, $M$ in \eqref{eq: relation
        vec(T)} is full-rank if and only if $Q$ in \eqref{eq: Q matrix} is full-rank. Next, we provide necessary and sufficient condition such that $Q$ is full-rank. Each block row of $Q$ can be written as
\begin{align}\label{eq: block row Q}
\begin{split}
 \xi(t)^{\T}\otimes C(t)&\!=\!
\begin{bmatrix}
 \xi_1(t)C(t) &  \xi_2(t)C(t)& \cdots & \xi_p(t)C(t)
\end{bmatrix}\\
&\!=\!
\underbrace{\begin{bmatrix}
 C(t) & \cdots & C(t)
\end{bmatrix}}_{F(t)}
\!\!\underbrace{\begin{bmatrix}
 \xi_1(t) I_p&  \cdots & 0 \\
 0         &  \cdots & 0\\
 \vdots & \ddots & \vdots\\
 0         & \cdots & \xi_p(t) I_p
\end{bmatrix}}_{\Xi(t)},
\end{split}
\end{align}
where $\xi_1,\cdots,\! \xi_p$ are the components of $\xi$ for $t\!\in\! \{N_0\!+\!1,\!\cdots , N\}$. Let $\Lambda=D\otimes I_p$. By noting that $\Xi(t+1)=\Lambda \Xi(t)$ and using the representation in \eqref{eq: block row Q}, we can write $Q$ in \eqref{eq: Q matrix}~as
\begin{align}\label{eq: W matrix}
 Q=
 \underbrace{\begin{bmatrix}
 F(N_0+1)\\
 F(N_0+2) \Lambda \\
 \vdots\\
 F(N) \Lambda^{N-N_0-1}
\end{bmatrix}}_{W}
\Xi(N_0+1),
\end{align}
From Assumption (A2), $\Xi(N_0+1)$ is full-rank. Hence, $Q$ in \eqref{eq: Q matrix} is full-rank if and only if the matrix $W$ in \eqref{eq: W matrix} is full-rank. Therefore, $M$ is full-rank if and only if $W$ is full-rank.
\end{IEEEproof}

\section{Numerical examples and comparisons}\label{sec:
  simulation}
We now present a series of numerical studies to validate the
effectiveness of the proposed optimization algorithm.

\subsection{Time-varying quadratic optimization}\label{subsec:
  time-varying quadratic term}
We consider a quadratic
time-varying optimization problem as $f(x)=x^\T H(t)x+b(t)^\T x$,
where $x\in\mathbb{R}^2$ and parameters $H(t)$ and $b(t)$ unknown. We
assume that $H(t) \succ 0$ for all $t\geq 0$ is a symmetric matrix. We express this
quadratic cost in the form \eqref{eq: optimization prob_intro} as
\begin{align}\label{eq:quad_cost}
f(x,z(t))=\begin{bmatrix} x_1 & x_2 & x_1^2 & 2x_1x_2 & x_2^2 \end{bmatrix} z(t),% ,
\end{align}
where $z(t)=[b_1(t), b_2(t), h_{11}(t), h_{12}(t), h_{22}(t)]^\T$.
The unknown parameter $z(t)$ obeys the dynamics in \eqref{eq:
  parameter dynamics block} with $A=\blkdiag([E,F])$, $E\!=\!
\begin{bmatrix} 0 & 1\\ -1& 0 \end{bmatrix}$, and
$F\!=\!\diag[0.98, 0.95, 0.981]$. We let
$z(0)\!=\![ -85.8, -77.9, 1047 , 329, 669]^{\T}$. The gradient of
\eqref{eq:quad_cost}~is
\begin{align}
  y(x,t)=
  \begin{bmatrix}
    1 & 0 & 2x_1 & 2x_2 & 0\\ 0 & 1 & 0 & 2x_1 & 2x_2
  \end{bmatrix}z(t).
\end{align}
We use $N_0=8$ and $N=26$ for the data collection as in~\eqref{eq:
  data}.
  
  \begin{figure}[!t]
  \centering
  \includegraphics[width=1\columnwidth,trim={0cm 0cm 0cm
    0cm},clip]{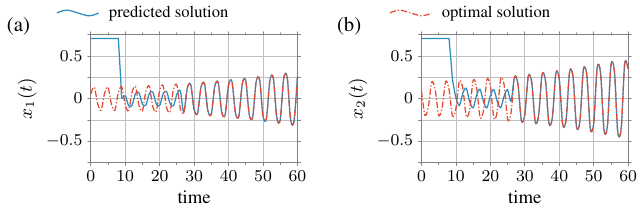}
  \caption{This figure shows the predicted (solid blue line) and the
    optimal solution (dashed red line) considering a quadratic cost
    function for the setting described in Section~\ref{subsec:
    time-varying quadratic term} as a function of time. Panels (a) and (b) 
    shows the evolution of the first and
    the second component of the predicted and the optimal solution,
    respectively. The data $X_0$ and $Y_0$ in \eqref{eq: data} are
    collected when $t\!\in\![0,8]$, where $[x_1(t),x_2(t)]^{\T}=
    [\sqrt{2}/2, \sqrt{2}/2]^{\T}$ and the data $X_1$ and $Y_1$ in 
    \eqref{eq: data} are collected when $t\! \in\! [9, 26]$, where $x(t)$ 
    is obtained using static gradient descent \eqref{eq: for comparison}. 
    For $t>26$, we predict the optimal solution using \eqref{eq:quad_sol}.
    We observe that the predicted value converges to the optimal solution 
    for $t>26$.}
    \label{fig: example_1}
\end{figure}

We estimate the realization of the system dynamics $A$ (using
$x(t)=[\sqrt{2}/2, \sqrt{2}/2]^{\T}$ for $t\in [0, N_0]$), up to
similarity transformation, using the collected data structured in the
Hankel matrix form as in \eqref{eq: hankel matrix}. Specifically, we
compute $\bar A$ using \eqref{eq: A_bar id} and then propagate
$\bar A$ to obtain
$[\bar{z}(N_0+1),\cdots, \bar{z}(N)]=\bar{A}^{p} [ \bar{z}(p-1),\cdots
,\bar{A}^{N-2p+1} \bar{z}(p-1)]$, where $\bar{z}(p-1)$ is computed
from \eqref{eq: Z_bar}. Then, from Theorems \ref{thrm: rank M} and
\ref{thrm: sufficient cond}, we have a full rank matrix $M$, which is
constructed as in \eqref{eq: relation vec(T)}. Using \eqref{eq:
  transformation matrix}, we can compute the transformation matrix $T$
that transforms the realization $\bar{A}$ to $A$ and $\bar{z}(t)$ to
the coordinate of $z(t)$ for $t\geq 0$.

Using the estimated $A$, we predict the evolution of $z(t)$ and
compute a solution to the optimization problem as
\begin{align}\label{eq:quad_sol}
  x(t)=-
  \underbrace{
  \begin{bmatrix}
    h_{11}(t) & h_{12}(t) \\ h_{12}(t) &  h_{22}(t)
  \end{bmatrix}^{-1}}_{H(t)^{-1}}
                                         \begin{bmatrix} b_1(t)\\
                                           b_2(t) \end{bmatrix}
  .
\end{align}
Fig.~\ref{fig: example_1} shows the predicted and the optimal solution
over time for the quadratic optimization. For comparison, we use a
static gradient descent algorithm that does not take into account the
time-varying nature of the cost function and updates the solution to
the optimization problem as
 \begin{align}\label{eq: for comparison}
 x(t)=x(t-1) - \eta y(x(t-1),t-1),
 \end{align}
 for $t\in [N_0+1, N]$, with $\eta=10^{-3}$, with
 $x(0)=[\sqrt{2}/2, \sqrt{2}/2]^{\T}$. As can be seen in Fig. \ref{fig:
   error_quad_poly}(a), our methods converges to the optimal solution
 in finite time, whereas the static gradient descent algorithm does
 not converge to the optimal solution.

\begin{figure}[!t]
  \centering
  \includegraphics[width=1\columnwidth,trim={0cm 0cm 0cm
    0cm},clip]{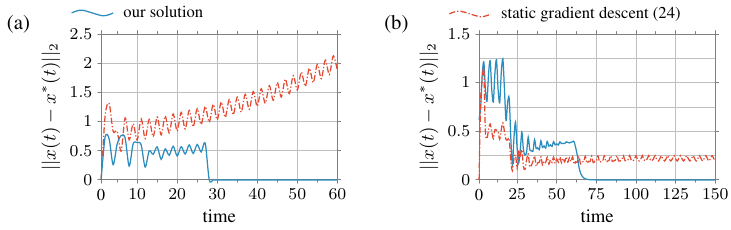}
  \caption{The figure shows a comparison of our proposed solution
    (solid blue line) with the static gradient descent algorithm \eqref{eq:
      for comparison} (dashed red line). Panel (a) and (b) shows the
    error in estimating the optimal solution for a quadratic cost
    function and a higher order polynomial, respectively. We observe
    that the error decays to zero for $t>N$ using our solution proposed in Section~\ref{subsec:
    time-varying quadratic term} and \ref{subsec:
    time-varying polynomial term}. However,
    for the static gradient descent \eqref{eq: for comparison}, we observe that the error diverges.}
    \label{fig: error_quad_poly}
\end{figure}

\subsection{Time-varying polynomial optimization}\label{subsec:
    time-varying polynomial term}
We consider now a third-order polynomial with time-varying
coefficients as
\begin{align}\label{eq:tensor_eq}
  f(x,z(t))=\sum_{r=1}^3 \underbrace{ \begin{bmatrix}  x^\T \otimes \cdots \otimes x^\T
    \end{bmatrix} }_{\text{r-terms}} z(t) ,
\end{align}
where $x\in \mathbb{R}^2$, and
$z(t)=[z_1(t)^\T , z_2(t)^\T , z_3(t)^\T]^\T$, with
$z_r(t)\in\mathbb{R}^{2^r}$ for $r\in\{1,2,3\}$ representing the
vector of coefficients. The dynamics of $z(t)$ follow \eqref{eq:
  parameter dynamics block} with $A:=\blkdiag([M_1,M_2,M_3])$,
$M_1=\begin{bmatrix} 0 & 1\\ -1 & 0 \end{bmatrix}$,
$M_2=\diag[0.98, 0.99, 0.99 , 0.95]$ and
$M_3=\diag[0.88 , 0.87 , 0.87 , 0.89 , 0.87 , 0.89 , 0.89 , 0.85]$.
We start from $z_1(0)= [-63.7, 110.2]^\T$, $z_2(0)=
\begin{bmatrix} 2.23, 2.46, 2.46, 6.24 \end{bmatrix}^\T$ and
$z_3(0)=[0.5, 0.3, 0.3, 0.4, 0.3, 0.4, 0.4, 0.6]^\T$.  The
gradient of \eqref{eq:tensor_eq} is
\begin{align*}
y(x,t)=\begin{bmatrix} P_1(x) & P_2(x) & P_3(x) \end{bmatrix} z(t),
\end{align*}
where $P_1(x)=I_2, P_2(x)=(x^\T\otimes I_2+I_2\otimes x^\T)$, and
$P_3(x)=(x^\T\otimes x^\T\otimes I_2+x^\T\otimes I_2\otimes x^\T+
I_2\otimes x^\T\otimes x^\T)$. We use $N_0=18$ and $N=60$ for the data
collection as in~\eqref{eq: data}.

\begin{figure}[!t]
  \centering
  \includegraphics[width=1\columnwidth,trim={0cm 0cm 0cm
    0cm},clip]{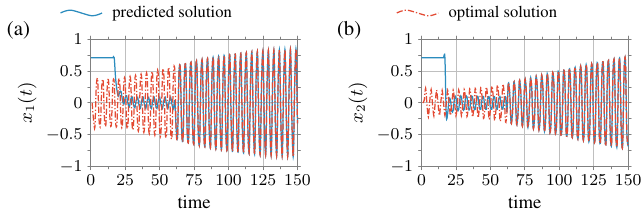}
  \caption{This figure shows the predicted (solid blue line) and the
    optimal solution (dashed red line) considering a third-order
    polynomial cost function for the setting described in Section~\ref{subsec:
    time-varying polynomial term} as a function of time. Panels (a) and (b) shows the evolution
    of the first and the second component of the predicted and the
    optimal solution, respectively.The data $X_0$ and $Y_0$ in \eqref{eq: data} are collected when $t\!\in\![0,18]$, where $[x_1(t),x_2(t)]^{\T}=[\sqrt{2}/2, \sqrt{2}/2]^{\T}$ and the data $X_1$ and $Y_1$ in \eqref{eq: data} are collected when $t\! \in\! [19, 60]$, where $x(t)$ is obtained using \eqref{eq: for comparison}. For $t>26$, we predict the optimal solution using algorithm~\ref{Algo1}. We observe that the predicted value converges to the optimal solution for $t>60$. }
    \label{fig: example_2}
\end{figure}

Unlike the case of quadratic cost function, a closed form expression
of the solution to the optimization problem is not available for
higher order polynomials. Instead, we employ a time-varying gradient
descent algorithm as in Algorithm~\ref{Algo1}, updating the solution
based on the parameter $z(t)$. We set the algorithm parameters as
$\beta=10^{-2}$, $D=500$, and $T=150$.  Fig.~\ref{fig: example_2}
shows the predicted and the optimal solution for a higher order
polynomial. We observe in Fig.~\ref{fig: example_2} that the predicted
solution attempts to track the optimal solution for $t\in [N_0+1, N]$
but never converges to it when computed via \eqref{eq: for
  comparison}. We also observe that the predicted value converges to
the optimal solution for $t >N$. Finally, Fig.~\ref{fig:
  error_quad_poly}(b) shows a comparison of the proposed algorithm and
the gradient descent algorithm using \eqref{eq: for comparison}. As
expected, our method converges to the optimal solution (up to the
convergence of the gradient descent algorithm), whereas the static
optimization algorithm \eqref{eq: for comparison} does not converge.

\begin{algorithm}[t!]
\SetKwInput{KwData}{Require}
  \label{Algo1}
  \SetAlgoLined
  \caption{\textbf{Time-varying~gradient descent} \cite{HI-TK-SI-AT:24}}
  \KwData{initial condition $x(N)$, step size $\beta$}
  \For{$t=N,N+1,\ldots,T$}{
Obtain $z(t)$ using $A$ and $z(0)$\\
Initialize $x_0(t)=x(t)$\\
\For{$d=0,1,\ldots,D-1$}{
\!\!$x_{d+1}(t) \! =\! x_d(t)\! -\! \beta \nabla_x f(x_d(t);z(t))$
}
Set $x(t+1)=x_D(t)$
}
\end{algorithm}

\subsection{Non-polynomial function}\label{section:non-polynomial}
\begin{figure}[!t]
  \centering
  \includegraphics[width=1\columnwidth,trim={0cm 0cm 0cm
    0cm},clip]{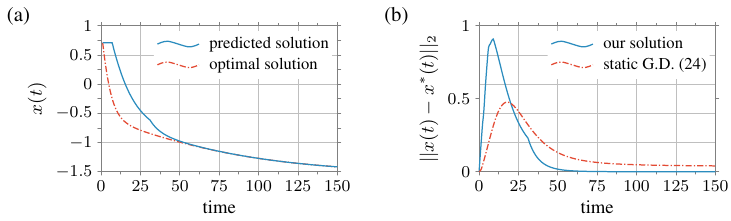}
  \caption{The figure in panel (a) shows the predicted (solid blue
    line) and the optimal solution (dashed red line) considering a
    non-polynomial cost function for the setting described in 
    Section~\ref{section:non-polynomial} as a function of time. 
    The data $X_0$ and $Y_0$ in \eqref{eq: data} are collected 
    when $t\!\in\![0,6]$, where $[x_1(t),x_2(t)]^{\T}=[\sqrt{2}/2, \sqrt{2}/2]^{\T}$ and the data $X_1$ and $Y_1$ in \eqref{eq: data} are collected when $t\! \in\! [7, 30]$, where $x(t)$ is obtained using \eqref{eq: for comparison}. For $t>26$, we predict the optimal solution using algorithm~\ref{Algo1}.
    The figure in panel (b) shows the error in
    computing the optimal solution, using our solution (solid blue line) and static gradient descent
    \eqref{eq: for comparison} (dashed red line). We observe that unlike the static gradient
    descent algorithm \eqref{eq: for comparison}, the predicted
    value converges to the optimal solution
    using our proposed solution.}
    \label{fig: example_3}
\end{figure}

In this section, we discuss the application of the proposed algorithm
to a non-polynomial problem. In particular, we consider the following
cost function:
\begin{align}\label{eq: non-poly}
f(x,z(t))=\begin{bmatrix} 2e^x & \sin x & x\end{bmatrix}z(t),
\end{align}
where $\map{f}{\mathbb{R} \times \mathbb{R}^3}{\mathbb{R}}$, and
$z(t)=[\theta_1(t), \theta_2(t), \theta_3(t)]^\T$ is a time varying
parameter. The dynamics of $z(t)$ is given by \eqref{eq: parameter
  dynamics block} with 
$A=\diag([0.99,0.97,0.98])$. The gradient of \eqref{eq: non-poly} is
\begin{align*}
y(x,t)=\begin{bmatrix} 2e^x & \cos x & 1\end{bmatrix}z(t).
\end{align*}
We collect the data as in \eqref{eq: data}, with $x(t)=x_0=0.7$ for
$t\in [0,6]$ and for $t\in[7,30]$ we compute $x(t)$ using \eqref{eq:
  for comparison} with $\eta=10^{-3}$.  Using the collected data, we
estimate the parameters $A$ and $z(t)$ for $t\geq 0$, following our
methodology as in the previous examples. Once $A$ is estimated, we
predict the evolution of $z(t)$ and compute a solution using
Algorithm~\ref{Algo1} with parameters $\beta=10^{-2}$, $D=500$, and
$T=150$. Fig. \ref{fig: example_3}(a) shows the predicted and the
optimal solution for the cost function in \eqref{eq: non-poly}
and Fig. \ref{fig: example_3}(b) compares the performance of
our solution using algorithm~\ref{Algo1} with the static gradient descent algorithm,
 as described in \eqref{eq: for comparison}.

\section{Conclusion}
In this work we study online time-varying optimization with
unknown, time-varying parameters, where gradient measurements are
leveraged to learn these parameters. We show how system identification
techniques can be used to estimate the parameters in their exact
coordinates through a transformation matrix. We also establish
necessary and sufficient conditions for computing the transformation
matrix, ensuring accurate parameter reconstruction. An open question
in our ongoing research is whether our approach can be extended to a
more general class of cost functions. For a general cost function, the
dynamics would be non-linear and time-varying. Exploring the
robustness of our method in noisy environments or including constrains
in our cost function are also important areas of future investigation.

\begin{appendix}
\setcounter{section}{5} 

\subsection{Technical lemmas}
\begin{appxlem}{\bf \emph{(Null-space of Kronecker product)}}\label{lemma: 1}
Let $A\in \mathbb{R}^{n\times m}$ and $B \in \mathbb{R}^{l\times k}$. Then,
\begin{align*}
\Null{A\otimes B} = \Null{A \otimes I_k} \cup \Null{I_m \otimes B}.
\end{align*}
\end{appxlem}
\begin{IEEEproof}
From \cite[Fact 7.4.23]{DSB:09}, we have
\begin{align*}
 \Range{A^{\T} \otimes B^{\T}} = \Range{A^{\T} \otimes I_k} \cap \Range{I_m \otimes B^{\T}}.
\end{align*}
The orthogonal complement of $\Range{A^{\T} \otimes B^{\T}}$ is written as
\begin{align*}
  \Range{A^{\T} \otimes B^{\T}}^{\perp} = \Range{A^{\T} \otimes I_k}^{\perp} \cup \Range{I_m \otimes B^{\T}}^{\perp}.
\end{align*}
The proof follows by noting that $\Range{X}^{\perp}=\Null{X^{\T}}$ for any matrix $X$.
\end{IEEEproof}
\begin{appxlem}{\bf \emph{(Basis of the null-space of Kronecker product)}}\label{lemma: 2}
Let $A\in \mathbb{R}^{n\times m}$ with $\Rank{\left(A\right)}=r$. Let $K \in \mathbb{R}^{m\times (m-r)}$ be a matrix whose columns span $\Null{A}$. Then, the columns of $K\otimes I_k$ span $\Null{A\otimes I_k}$.
\end{appxlem}
\begin{IEEEproof}
We begin  by noting that $\col{\left(K\otimes I_k\right)} \!\in\! \Null{A\otimes I_k}$ since $\left(A\otimes I_k\right)\left(K\otimes I_k\right)\! =\! AK \otimes I_k \!=\!0$. Thus, the proof boils down to showing that $\Rank{\left(K \otimes I_k\right)}\!=\!\dim{\left(\Null{A\otimes I_k}\right)}$. Using \cite[Fact 7.4.24]{DSB:09},~we have $\Rank{\left(A\otimes I_k \right)}\! =\! rk$. Consequently, using the rank-nullity theorem, we have, $\dim{\left(\Null{A\otimes I_k}\right)}\!=\!mk\!-\!rk$. Further, since $\Rank{\left(K\right)}\!=\!m-r$, using \cite[Fact 7.4.24]{DSB:09} we have, $\Rank{\left(K \otimes I_k\right)\!=\!mk-rk}\!=\!\dim{\left(\Null{A\otimes I_k}\right)}$. Therefore, the columns of $K \otimes I_k$ span $\Null{A\otimes I_k}$.
\end{IEEEproof}

\end{appendix}

%\clearpage
\bibliographystyle{unsrt}
\bibliography{alias,Main,FP,New}

\end{document}